\documentclass[12pt]{article}

\usepackage{epsfig,cite}
\usepackage{graphics}
\usepackage{amssymb,amsmath}
\usepackage{times,comment}
\newcommand{\eqna}[1]{\begin{subequations} \label{#1}
\begin{eqnarray}}
\def\eena{\end{eqnarray}
\end{subequations}}

\def\be{\begin{equation}}
\def\eqn#1{\be\label{#1}}
\def\ee{\end{equation}}

\def\bea{\begin{eqnarray}}
\def\eqnn#1{\bea\label{#1}}
\def\eea{\end{eqnarray}}

\def\nn{\nonumber}

\def\nd{\end{document}}

\input epsf.tex
\newcount\figno
\def\fig#1#2#3{
\par\begingroup\parindent=0pt\leftskip=1cm\rightskip=1cm\parindent=0pt
\baselineskip=11pt \global\advance\figno by 1 
\epsfxsize=#3 \centerline{\epsfbox{#2}} \vskip 12pt
#1\par
\endgroup\par}
\def\figlabel#1{\xdef#1{\the\figno}}
\def\encadremath#1{\vbox{\hrule\hbox{\vrule\kern8pt\vbox{\kern8pt
\hbox{$\displaystyle #1$}\kern8pt} \kern8pt\vrule}\hrule}}

  \def\tV{{\tilde V}}

\def\bu{\noindent $\bullet~$}
\def\rank{{\rm rank}}

\def\riga{-\kern-4pt - \kern-4pt -}
\font\fat=cmsy10 scaled\magstep5

\def\Bbullet{\raise-3pt\hbox{\fat\char"0F}}

\def\mt{\mapsto}

\def\Box{
\vbox{ \halign to5pt{\strut##& \hfil ## \hfil \cr &$\kern -0.5pt
\sqcap$ \cr \noalign{\kern -5pt \hrule} }}~}

\def\down{\raise1.5pt\hbox{$\phantom{a}_2$}\downarrow}

\def\downa{\raise1.5pt\hbox{$\phantom{a}_{2\atop m_2}$}\downarrow}

\def\llr{\longrightarrow}

\def\({\left(}
\def\){\right)}
\def\eps{\epsilon}

\def\ha{{\textstyle{\frac{1}{2}}}}

\def\bbz{\mathbb{Z}}
\def\bbc{\mathbb{C}}
\def\bbr{\mathbb{R}}
\def\bbn{\mathbb{N}}
\def\a{\alpha}
\def\b{\beta}

\def\vr{\vert}

\def\l{\lambda}

\def\D{{\Delta}}

\def\ca{{\cal A}} \def\cb{{\cal B}} \def\cc{{\cal C}}
\def\cd{{\cal D}}  \def\cf{{\cal F}}
\def\cg{{\cal G}} \def\ch{{\cal H}} 
 \def\ck{{\cal K}} 
\def\cm{{\cal M}} \def\cn{{\cal N}} 
\def\cp{{\cal P}}  \def\car{{\cal R}}
 \def\ct{{\cal T}} \def\cu{{\cal U}}

\def\ido{intertwining differential operator}

\def\L{\Lambda}
\def\r{\rho}

\begin{document}
\sloppy \raggedbottom

\begin{center}

\textsf{\Large\bf Parabolic Verma Modules and \\[3pt]
Invariant Differential Operators
\footnote{Plenary talk at  Bogolyubov Conference, Moscow-Dubna, 11.9.2019.}}

 \vspace{7mm}

{\large\bf  V.K. Dobrev}

\vspace{3mm}

 \emph{Institute of Nuclear Research and Nuclear Energy,\\
 Bulgarian Academy of Sciences, \\
72 Tsarigradsko Chaussee, 1784 Sofia, Bulgaria}

\end{center}

\begin{abstract}
In the present paper we continue the project of systematic classification and
construction of invariant differential operators for non-compact semisimple Lie groups. This time we make the stress on one of the main
building blocks, namely the Verma modules and the corresponding parabolic subalgebras.
In particular, we start the study of the relation between the parabolic subalgebras of real semisimple Lie algebras
and of their complexification. Two cases are given in more detail: the conformal algebra of 4D Minkowski space-time and the minimal parabolics of classical real semisimple Lie algebras.
 \end{abstract}

\section{Introduction}

Verma modules play crucial role in the  classification and
construction of invariant differential operators for non-compact semisimple Lie groups.
In our approach we make important use of the relations between the parabolic subalgebras of
a real semisimple Lie algebra $\cg_0$ and those of its complexification $\cg$. It the present paper
we make this relation more transparent briefly in general and in some detail in examples.

\section{Preliminaries}

 Let $G$ be a semisimple non-compact Lie group, and $K$ a
maximal compact subgroup of $G$. Then we have the {\it Iwasawa
decomposition} ~$G=KA_0N_0$, where ~$A_0$~ is Abelian simply
connected vector subgroup of ~$G$, ~$N_0$~ is a nilpotent simply
connected subgroup of ~$G$~ preserved by the action of ~$A_0$.
Further, let $M_0$ be the centralizer of $A_0$ in $K$. Then the
subgroup ~$P_0 ~=~ M_0 A_0 N_0$~ is a {\it minimal parabolic subgroup} of
$G$.  A {\it parabolic subgroup} ~$P ~=~ M A N$~ is any subgroup of $G$
which contains a minimal parabolic subgroup. In the general case ~$A \subset A_0$~ is abelian,
~$N\subset N_0$~ is a nilpotent simply
connected subgroup of ~$G$~ preserved by the action of ~$A$, $M\supset M_0$ is a maximal centralizer of $A$ in $G$.

Further, let ~$\cg_0,\ck,\cp,\cm,\ca,\cn$~ denote the Lie algebras of ~$G,K,P,M,A,N$, resp.

We note also another extremal case :  ~{\it maximal   parabolic
subgroup}~  when  $\rank\, A =1$,   ~{\it maximal   parabolic
subalgebra}  when ~$\dim\, \ca=1$.

Let ~$\nu$~ be a (non-unitary) character of ~$A$, ~$\nu\in\ca^*$.
 Let ~ $\mu$ ~ fix a finite-dimensional (non-unitary) representation
~$D^\mu$~ of $M$ on the space ~$V_\mu\,$.
 In the case when $M$ is cuspidal then we may  use also  the
 discrete series representation of $M$ with the same Casimirs as $D^\mu$.

 We call the induced
representation ~$\chi =$ Ind$^G_{P}(\mu\otimes\nu \otimes 1)$~ an
~{\it \it elementary representation} of $G$ \cite{DMPPT}. (These are
called {\it generalized principal series representations} (or {\it
limits thereof}) in \cite{Knapp}.)   Their spaces of functions are:  \eqn{func}
\cc_\chi ~=~ \{ \cf \in C^\infty(G,V_\mu) ~ \vr ~ \cf (gman) ~=~
e^{-\nu(H)} \cdot D^\mu(m^{-1})\, \cf (g) \} \ee where ~$a=
\exp(H)\in A$, ~$H\in\ca\,$, ~$m\in M$, ~$n\in N$. The
representation action is the {\it left regular action}:  \eqn{lrega}
(\ct^\chi(g)\cf) (g') ~=~ \cf (g^{-1}g') ~, \quad g,g'\in G\ .\ee

\section{Verma modules}

 An important ingredient in our considerations are the ~{\it
highest/lowest weight representations}~ of ~$\cg$, where ~$\cg$~ is the complexification of ~$\cg_0$. These can be
realized as (factor-modules of) Verma modules ~$V^\L$~ over
~$\cg$, where ~$\L\in (\ch)^*$, ~$\ch$ is a Cartan
subalgebra of ~$\cg$, the weight ~$\L = \L(\chi)$~ is determined
uniquely from $\chi$ \cite{Dob}.

We recall that a {\it {Verma module}}  \cite{Dix,Ver} is a  $\cg$--module induced from a character of a Borel subalgebra ~$\cb=\ch\oplus \cn$.
(We use the standard triangular decomposition ~$\cg ~= \cn\oplus \ch \oplus \cn^-$.) Let $1_\L$ denote
the $\cb$--module $1_\L = \bbc v_0$ such that $\cn v_0 = 0$, $H
v_0 = \L(H) v_0$, ~$H\in\ch$, $(\dim 1_\L = 1)$. Let $V^\L$ be the corresponding
{Verma} module, then \eqn{verl} V^\L = U(\cg)\otimes_{U(\cb)}1_\L
~.\ee Using \eqn{decucgb}  U(\cg) = U(\cn^-)\otimes U(\cb) \ee one
has \eqn{verla} V^\L = U(\cn^-)\otimes 1_\L ~. \ee Obviously $V^\L$
is a highest weight module with highest weight $\L$, and highest
weight vector $v_0$, $V^\L\cong U(\cn^-)$ as vector spaces.

Actually, since our ERs  are induced from finite-dimensional
representations of ~$\cm$~ (or their limits) the Verma modules are
always reducible. Thus, it is more convenient to use ~{\it
generalised Verma modules} ~$\tV^\L$~ such that the role of the
highest/lowest weight vector $v_0$ is taken by the
(finite-dimensional) space ~$V_\mu\,v_0\,$.

It is important to note that the generalized Verma modules defined just above and related to ERs are special cases
of ~{\it parabolic Verma modules (PVM)}~ which are introduced in purely algebraic context.
More precisely, their construction is as follows. (Below we define PVM using \cite{Lep}, except that there these modules are called generalized Verma modules (which in our approach  is used differently, see above and \cite{VKD1}).

Let ~$\D$~ be the root system of ~$(\cg,\ch)$, ~$\D_+,\D_-$~ denote the positive, negative, roots. Let ~$\a_i\in\ch^*$~ be the simple roots where ~$i=1,\ldots,\ell ~=~ \dim\ch$.  Let ~$e_i$ (resp. $f_i$) be a non-zero element of the root space ~$\cg_{\a_i}$ (resp. ~$\cg_{-\a_i}$) for all $i=1,\ldots,\ell$ normalized so that ~$[e_i,f_i]=h_i$, and ~$\a_i(h_i)=2$.
Let ~$S$~ be a subset of the set ~$J=\{1,\ldots,\ell\}$.
Let ~$\cg_S$~ be the subalgebra of $\cg$ generated by  ~$\{h_i,e_i,f_i\}_{i\in S}\,$;  ~$\ch_S$~  the span of
~$\{h_i\}_{i\in S}\,$; ~$\D^S = \D\cap \amalg_{i\in S} \bbz\a_i\,$; ~$\D^S_\pm = \D_\pm \cap \D^S$;
~$\D(S)_\pm = \D_\pm -  \D^S_\pm\,$.

Further, we define  the following subalgebras  of $\cg$~:~ $\cn = \amalg_{\b\in\D_+}\cg^\b$; ~
$\cn^- = \amalg_{\b\in\D_-}\cg^\b$ (the latter two were already used above); ~$\cn_S = \amalg_{\b\in\D_+^S}\cg^\b$;
 ~$\cn_S^- = \amalg_{\b\in\D_-^S}\cg^\b$; ~$\cu_S = \amalg_{\b\in\D_+(S)}\cg^\b$; ~
 ~$\cu^-_S = \amalg_{\b\in\D_-(S)}\cg^\b$; ~$\car_S = \cg_S \oplus \ch$;
 ~$\cp_S = \car_S \oplus \cu_S\,$.

 Then we have: ~$\cg ~=~ \cn \oplus \ch \oplus\cn^-$, ~$\cg_S ~=~ \cn_S \oplus \ch_S \oplus\cn^-_S$;
 ~$\cn = \cn_S \oplus \cu_S\,$; ~$\cn^- = \cn^-_S \oplus \cu^-_S\,$;~ $\car_S = \cn_S \oplus \ch \oplus \cn^-_S\,$;
 ~$\cg ~=~ \cu^-_S \oplus \cp_S\,$. \\
 Further, ~$\cg_S$~ is a split semisimple Lie algebra with splitting Cartan subalgebra ~$\ch_S\,$; ~$\car_S$~ is a reductive Lie algebra with commutator subalgebra ~$\cg_S$~ and centre a subalgebra of ~$\ch$.

 As ~$S$~ varies among the subsets of ~$J$, ~$\cp_S$~ varies among the parabolic subalgebras of ~$\cg$~ containing the Borel subalgebra ~$\cb = \ch \oplus \cn$. The reductive part of ~$\cp_S$~ is ~$\car_S$~ and the nilpotent part of
 ~$\cp_S$~ is ~$\cu_S\,$. We note that if ~$S=\emptyset$, then ~$\cp_\emptyset ~=~ \cb$ (since $\cg_\emptyset ~=~ 0$,
 $\car_\emptyset ~=~ \ch$, ~ $\cu_\emptyset ~=~ \cn$).

 Now let ~$P$~ indexes the set of (equivalent classes of) finite-dimensional irreducible ~$\cg$-modules in the usual way - via the highest weight.  Let ~$P_S ~=~ \{\l\in\ch^* | \l(h_i) \in \bbz_+ ~{\rm for~all}~ i\in S\}$.
 Then it is clear that there is a natural bijection,  ~$\l \mt M(\l)$, between ~$P_S$~ and the set of
 (equivalent classes of) finite-dimensional irreducible ~$\car_S$-modules which are irreducible as ~$\cg_S$-modules.

  For any ~$\l\in P_S$, denote by ~$V_S^{M(\l)}$~ the corresponding ~{\it parabolic Verma module} (PVM)~:~ the ~$\cg$-module
 ~ind$(M(\l),\cg)$~ induced by the ~$P_S$-module ~$M(\l)$, viewed as an ~$\car_S$-module in the natural way and as a trivial ~$\cu_S$-module.

 Note that for ~$S=\emptyset$~ the parabolic Verma modules coincide with the usual  Verma modules~:~
 $V_\emptyset^{M(\l)} ~=~ V^\l$.

Further we discuss the reducibility of  Verma modules.

A classic result of \cite{BGG} states that a Verma module ~$V^\L$~ is reducible iff
\eqn{bggr} (\L+\r, \b^\vee ) ~=~ m \ , \quad m\in\bbn, ~~\b\in\D^+, ~~
 \ \ee
where ~$\b^\vee \equiv 2 \b /(\b,\b)$, ~$\r$~ is half the sum of the positive roots of
~$\cg$.

The same criterion of reducibility is valid for  generalized Verma modules, though it is trivially satisfied
for the ~$\cm$-compact roots, and is essential only for ~ $\cm$-non-compact roots. (We recall that
~$\cm$-compact roots are those elements of $\D$ that belong to the root system of $\cm^\bbc$, the latter being identified as a subset of $\D$.)

The same criterion of reducibility is valid for a parabolic Verma module ~$V_S^{M(\L)}$~ though it is trivially satisfied for ~$\b \in \D_+^S$, and is essential only for  ~$\b\in \D_+(S)$.

 When \eqref{bggr}  holds then the Verma module with shifted
weight ~$V^{\L-m\b}$ (or ~$\tV^{\L-m\b}$ ~ for GVM and $\b$
$\cm$-non-compact, or ~$V_S^{M(\L-m\b)}$~ for PVM and ~$\b\in \D_+(S)$) is embedded in the Verma module ~$V^{\L}$ (or
~$\tV^{\L}$, or ~$V_S^{M(\L)}$).

The above embedding is realized by a singular vector
~$v_s$~ determined by a polynomial ~$\cp_{m,\b}(\cn^-)$~ in the
universal enveloping algebra ~$(U(\cn^-))\ v_0\,$.
 More explicitly, \cite{Dob}, ~$v^s_{m,\b} = \cp_{m,\b}\, v_0\,$.
 Relatedly, then
there exists \cite{Dob} an {\it \ido}  \eqn{invop}  \cd_{m,\b} ~:~ \cc_{\chi(\L)}
~\llr ~ \cc_{\chi(\L-m\b)} \ee given explicitly by: \eqn{singvv}
 \cd_{m,\b} ~=~ \cp_{m,\b}(\widehat{\cn^-})  \ee where
~$\widehat{\cn^-}$~ denotes the {\it right action} on the functions
~$\cf$.

 In the next Section we shall consider the example of the conformal algebra.

 \section{Conformal algebra}

 The conformal algebra in four-dimensional space-time  is ~$\cg_0 = su(2,2)$ ($\cong so(4,2)$).
  It has three nonconjugate parabolic
subalgebras ($\cp = \cm \oplus \ca \oplus \cn$):
\eqnn{parab}
\cp_{0} ~&=&~ so(2) \oplus \ca_0 \oplus \cn_0\ , \nn\\
&& \dim\ca_0 =2\ ,\quad \dim\cn_0 = 6\ , \nn\\
\cp_{1} ~&=&~ so(2) \oplus sl(2,\bbr) \oplus \ca_1 \oplus \cn_1\ , \\
&& \dim\ca_1 =1\ ,\quad \dim\cn_1 = 5\ , \nn\\
\cp_{2} ~&=&~ so(3,1) \oplus \ca_2 \oplus \cn_2\ , \nn\\
&& \dim\ca_2 =1\ ,\quad \dim\cn_2 = 4\ , \nn\eea
where ~$\cp_0$~ is the minimal parabolic, ~$\cp_1$~ is maximal cuspidal,
~$\cp_2$~ is maximal noncuspidal.

  \subsection{Maximal non-cuspidal case}

 We consider
  the following   Bruhat decomposition \cite{Bru} (consistent with the maximal non-cuspidal parabolic subalgebra $\cp_2$):
 \eqn{bruh} \cg_0 ~=~ \cg^+_0 \oplus \cm_2 \oplus \ca_2 \oplus \cg^-_0  \ee
 where ~$\cm_2$~ is the six-dimensional Lorentz subalgebra $so(3,1)$, ~$\ca_2$~ is the dilatation subalgebra,
 ~$\cg^+_0,\cg^-_0$~ is the four-dimensional isomorphic translation subalgebra, resp., special conformal transformations subalgebra.

  In this case the ERs of ~$su(2,2)$~are parametrized by triples: ~$\chi = [j_1,j_2;d]$, where ~$j_1,j_2\in \ha\bbz_+$~ parametrize
 the finite-dimensional representations of   ~$\cm_2$, while the number ~$d$~
 parametrizing the representations of   ~$\ca_2$~  is  called the {\it conformal weight} or
energy.

 The  complexificaton of ~$\cg_0$~ is ~$\cg = sl(4)$. The root system of ~$\cg$~ is given by:
 \eqn{roots}
 \D^+ = \{\a_1,\a_2,\a_3,\a_{12}= \a_1+\a_2,  \a_{12}= \a_1+\a_2, \a_{23}= \a_2+\a_3, \a_{13}= \a_1+\a_2+\a_3 \} \ee
 where ~$\{\a_1,\a_2,\a_3\}$~ are the simple roots.

 Note that when relating the root systems of ~$\cg$~ to ~$\cg_0$~ relative to the Bruhat decomposition \eqref{bruh} the roots ~$\{\a_1,\a_3\}$~ are  $\cm_2$-compact, the rest
 are $\cm_2$-non-compact.

The reducibility conditions \eqref{bggr} of a Verma module ~$V^\L$~ over ~$\cg$~ are written explicitly as follows:
\eqna{redu} m_1 ~&=&~ (\L+\r, \a_1 ) \in \bbn \\
m_2 ~&=&~ (\L+\r, \a_2 ) \in \bbn \\
m_3 ~&=&~ (\L+\r, \a_3 ) \in \bbn \\
m_{12} ~&=&~ (\L+\r, \a_{12} ) = m_1+m_2 \in \bbn \\
m_{23} ~&=&~ (\L+\r, \a_{23} ) =m_2+m_3\in \bbn \\
m_{13} ~&=&~ (\L+\r, \a_{13} ) = m_1+m_2+m_3 \in \bbn \eena

We want to apply these conditions to the signatures $\chi$ of the ERs.
In these terms we have \cite{VKD1}:
 \eqna{reduc} m_1 ~&=&~  2j_1+1 \in \bbn \\
m_2 ~&=&~ 2- d-j_1-j_2  \in \bbn \\
m_3 ~&=&~ 2j_2+1  \in \bbn \\
m_{12} ~&=&~  3- d+j_1-j_2    \in \bbn \\
m_{23} ~&=&~ 3- d-j_1+j_2  \in \bbn \\
m_{13} ~&=&~ 4- d+j_1+j_2   \in \bbn \eena

Note that (\ref{redu}a,c) are fulfilled always since ~$2j_1+1\in\bbn$,
~$2j_2+1\in\bbn$, as expected for the $\cm_2$-compact roots.
On the other hand the expressions in the other cases depend on $d$ and may be arbitrary.

Note that ~$m_i$~ considered abstractly are called Dynkin labels, while together with ~$m_{ij}$~
they are called Harish-Chandra parameters \cite{Har}:
\eqn{dynhc} m_\b \equiv (\L+\r, \b )\ ,  \ee where $\b$ is any
positive root of $\cg$. These parameters are redundant, since
they are expressed in terms of the Dynkin labels, however,   some
statements are best formulated in their terms.

\subsection{Maximal cuspidal case}

We consider again  the conformal algebra  ~$\cg_0 = su(2,2)$ ($\cong so(4,2)$). Here we consider
  the following   Bruhat decomposition (consistent with the maximal cuspidal parabolic subalgebra):
  \eqn{bruh1} \cg_0 ~=~ \cg^+_1 \oplus \cm_1 \oplus \ca_1 \oplus \cg^-_1  \ee
 where ~$\cm_1 = so(2) \oplus so(2,1)$, ~$\ca_1$~ is one-dimensional,
 ~$\cg^+_1,\cg^-_1$~ are  five-dimensional isomorphic   subalgebras.

 The signatures of the ERs in this case are \cite{VKD-cusp}:
\eqn{sign} \chi_1 ~=~ \{\,n',k,\eps,\nu'\,\}\ ,\ee
 where ~$n'\in\bbz$~ is a character of ~$so(2)$, ~$\nu'\in\bbc$~ is a
character of ~$A_1\,$, ~$k,\eps$~ fix a discrete series
representation of ~$so(2,1)$, ~$k\in\bbn$, ~$\eps =\pm1$, or a
limit thereof when $k=0$.

The relation with the ~$sl(4)$~ Dynkin labels is as follows \cite{VKD-cusp}:
\eqn{intp} m_1 ~=~ \ha (k-\nu'+n') \ ,\quad m_2 ~=~ -k
\ , \quad m_3 ~=~ \ha (k-\nu' -n')\ . \ee
For the analysis we need the additional Harish-Chandra parameters:
\eqn{harch} m_{12} ~=~ \ha (n'-k-\nu') , ~~ m_{23} ~=~ - \ha (k+\mu'+n') , ~~
m_{13} ~=~ -\nu' \ee

We see that if ~$\nu'\notin \bbz$~ then no Harish-Chandra parameter can be a positive integer, thus, the ERs would be irreducible.

Thus, we consider the case ~$\nu'\in \bbz$. Actually, we shall use the analysis
of the partially equivalent ERs in this case done in \cite{VKD-cusp}. Thus, we use a parametrization taken from there
 (up to change of sign) by three positive integers: ~$p,\nu,n$, so that  we have:
\eqnn{relatt} && \l_{p,\nu,n} ~=~ (m_1,m_2,m_3)_{p,\nu,n}  ~=~ (-p-\nu,  \nu , -n-\nu) , \\
&& m_{12} ~=~ -p , ~~ m_{23} ~=~ - n , ~~ m_{13} ~=~ -p-n-\nu \nn\eea
It is known that when relating the root systems of ~$\cg$~ to ~$\cg_0$~ relative to the Bruhat decomposition \eqref{bruh1} the root ~$\a_2$~ is   compact, the rest  are non-compact, and the above parametrization is
consistent with this.

\section{Parabolic Verma modules for $sl(4)$}

Here we enumerate the parabolic Verma modules for $sl(4)$. For this we need to produce the list of the various   ~$\D^S_\pm\,$ and the corresponding parabolic subalgebras. We have:
 \eqnn{setss} && \D^\emptyset_\pm ~=~ 0, ~~~ \D^1_\pm ~=~ \{\pm\a_1 \},  ~~~\D^2_\pm ~=~ \{\pm\a_2\},  ~~~\D^3_\pm ~=~ \{\pm\a_3\},
 ~~~\D^{123}_\pm ~=~ \D_\pm \ ,
 \nn\\
&& \D^{12}_\pm ~=~ \{\pm\a_1,\pm\a_{2},\pm\a_{12}\},
~~~ \D^{23}_\pm ~=~ \{\pm\a_2,\pm\a_{3}, \pm\a_{23}\}, ~~~\D^{1,3}_\pm ~=~ \{\pm\a_1, \pm\a_{3}\}
\eea
\eqnn{setsss} &&   P_\emptyset ~=~ \cb, ~~~   P_{\{i\}} ~=~ \cg^i_- \amalg \cb , ~~(i=1,2,3)\\
&&P_{\{12\}} ~=~ \cg^1_- \amalg \cg^2_- \amalg\cg^{12}_- \amalg \cb ,
~~~P_{\{23\}} ~=~ \cg^2_- \amalg \cg^3_- \amalg\cg^{23}_- \amalg \cb ,  \nn\\
&&P_{\{1,3\}} ~=~ \cg^1_- \amalg \cg^3_-  \amalg \cb ,
~~~P_{\{123\}} ~=~ \cg_- \amalg  \cb ~=~ \cg \nn\eea

Now we can make connection with some generalized Verma modules.

In order to compare the parabolic subalgebras of the real form with the parabolic subalgebras of the
complexification ~$sl(4,\bbc)$~ we need the complexification of \eqref{parab}. We have:
\eqna{parabc}
\cp^\bbc_{0} ~&=&~ so(2,\bbc) \oplus \ca^\bbc_0 \oplus \cn^\bbc_0    , \\
 \cp^\bbc_{1} ~&=&~ so(2,\bbc) \oplus sl(2,\bbc) \oplus \ca^\bbc_1 \oplus \cn^\bbc_1\ , \\
 \cp^\bbc_{2} ~&=&~ so(4,\bbc) \oplus \ca^\bbc_2 \oplus \cn^\bbc_2\ \eena

\bu First  we note that ~$\ch \cong so(2,\bbc) \oplus \ca^\bbc_0\,$, ~$\cn^\bbc_0 ~\cong~  \cn$, ~thus:
 ~$\cp^\bbc_{0}  ~=~ \cb$.

 \bu  Further, we record the triangular decomposition ~:~
   $sl(2,\bbc) ~=~ sl(2,\bbc)^+\oplus sl(2,\bbc)_h \oplus sl(2,\bbc)^-$,
   where  ~$sl(2,\bbc)_h$~ is a Cartan subalgebra of ~$sl(2,\bbc)$.
Then we    note that ~$\ch ~\cong ~so(2,\bbc) \oplus sl(2,\bbc)_h \oplus \ca^\bbc_1$,
~$\cn ~\cong~ sl(2,\bbc)^+ \oplus \cn^\bbc_1\,$,   and thus we have:
\eqn{parabd}  \cp^\bbc_{1} ~=~ sl(2,\bbc)^- \oplus \cb ~\cong~ P_{\{2\}}\ee
Now we can note that the GVM ~$V^{\L(\chi_1)}$  with ~$\chi_1 = \{\,n',k,\eps,\nu'\}$~ is isomorphic to PWM ~$V_{2}^{M(\l_{p,\nu,n})}$, so that  \eqn{lammu} \l (h_1,h_2,h_3) ~=~ (m_1-1,m_2-1,m_3-1) ~=~ (-p-\nu-1,\nu-1,-n-\nu-1) \ee

\bu Next, we record the triangular decomposition ~:~
   $so(4,\bbc) ~=~ so(4,\bbc)^+\oplus so(4,\bbc)_h \oplus so(4,\bbc)^-$,
   where  ~$so(4,\bbc)_h$~ is a Cartan subalgebra of ~$so(4,\bbc)$.
 Next  we note that ~$\ch ~\cong ~so(4,\bbc)_h  \oplus \ca^\bbc_2$.
  ~$\cn ~\cong~ so(4,\bbc)^+ \oplus \cn^\bbc_2\,$,   and thus we have:
\eqn{parabd}  \cp^\bbc_{2} ~=~ so(4,\bbc)^- \oplus \cb ~\cong~ P_{\{1,3\}}\ee
Finally,  we note that the GVM ~$V^{\L(\chi)}$  with ~$\chi = [j_1,j_2; d]$~ is isomorphic to PWM ~$V_{1,3}^{M(\l)}$, where \eqn{lamm} \l (h_1,h_2,h_3) ~=~ (m_1-1,m_2-1,m_3-1) ~=~ (2j_1, 1-d-j_1-j_2, 2j_2) \ee

\section{Minimal parabolics vs. complex parabolic subalgebras}

Here we briefly discuss the relation of minimal parabolics ~$\cp_0=\cm_0\oplus \ca_0 \oplus \cn_0$~  of ~{\it classical}~ real Lie algebras to the parabolic subalgebras of their complexification. (For the minimal parabolic subalgebras and the enumeration of simple roots we use info from \cite{VKD1,VKD-par}.)

In the case of split real Lie algebras ~$\cg_r$~ we have the general fact that ~$\cm_0=0$~  and then the complexification of the minimal parabolic of  ~$\cg_r$~ is isomorphic to the Borel subalgebra of ~$\cg_r^\bbc\,$.
(We recall also that in this case ~$\ca_0^\bbc \cong \ch(\cg_r^\bbc)$.)
 For completeness we list the
classical split real Lie algebras: ~$sl(n,\bbr)$, ~$so(r,r)$, ~$so(r+1,r)$,  ~$sp(n,\bbr)$.

There are cases of non-split  real Lie algebras ~$\cg_r$~ when the minimal parabolic is isomorphic to
~$\cb(\cg_r^\bbc)$. That is when the subalgebra ~$\cm_0$~ is abelian. Then
~$\cm_0^\bbc \oplus \ca_0^\bbc \cong \ch(\cg_r^\bbc)$. In the classical case this is the real Lie algebra
~$su(n,n)$, ($n>1$), where  we have: ~$\cm_0 ~=~ u(1) \oplus \cdots \oplus \cdots u(1)$, ($n-1$~ entries), ~$\dim_\bbr \ca_0 ~=~n$,
~$\dim_\bbr \cn_0 ~=~n(2n-1)$.  Then,
~$\cp_0^\bbc ~\cong ~\cb(sl(2n,\bbc))$.

Next we consider the rest of the real Lie algebras where the relation of the minimal parabolic to the
complex parabolics is more involved.

In the case of ~$su^*(2n)$~ ($n>1$) the minimal parabolic subalgebra   is given by: ~$\cm_0 ~=~ su(2) \oplus \cdots \oplus \cdots su(2)$, ($n$~ entries), ~$\dim_\bbr \ca_0 ~=~n-1$,
~$\dim_\bbr \cn_0 ~=~2n(n-1)$. Thus,
~$\cp^\bbc_0 \cong P_{1,3,\ldots,2n-1} \cong \amalg_{i=1}^{n} \cg_-^{2i-1}\amalg \cb $.

In the  case of ~$su(p,r)$~ ($p>r\geq1$) the minimal parabolic subalgebra   is given by:  ~$\cm_0 ~=~ su(p-r)\oplus u(1) \oplus \cdots \oplus \cdots u(1)$, ($r$~ entries), ~$\dim_\bbr \ca_0 ~=~r$, ~$\dim_\bbr \cn_0 ~=~r(2p-1)$. Thus,
~$\cp^\bbc_0 \cong P_{r+1,\ldots,p-1} \cong \amalg_{i=r+1}^{p-1} \cg_-^{i}\amalg \cb $.

In the case of ~$so(p,r)$~ ($p> r+1$) the minimal parabolic subalgebra   is given by:  ~$\cm_0 ~=~ so(p-r)$, ~$\dim_\bbr \ca_0 ~=~r$, ~$\dim_\bbr \cn_0 ~=~r(p-1)$.
Thus,
~$\cp^\bbc_0 \cong P_{r+1,\ldots,[(p+r)/2]} \cong \amalg_{i=r+1}^{[(p+r)/2]} \cg_-^{i}\amalg \cb $.

Next we consider ~$sp(p,r)$~ ($p\geq r$). The minimal parabolic subalgebra   is given by:  ~$\cm_0 ~=~ sp(p-r) \oplus sp(1) \oplus \cdots \oplus sp(1)$, ($r$ factors), ~$\dim_\bbr \ca_0 ~=~r$, ~$\dim_\bbr \cn_0 ~=~ r(4p-1)$.
In the case ~$p=r$~ we have ~$\cp^\bbc_0 \cong P_{1,3,\ldots,2r-1} \cong \amalg_{i=1}^r \cg_-^{2i-1}\amalg \cb $.
In the case ~$p>r$~ we have ~$\cp^\bbc_0 \cong P_{1,3,\ldots,2r-1,2r+1,2r+2,\ldots,p+r}
\cong \amalg_{i=1}^r \cg_-^{2i-1}\amalg_{j=2r+1}^{p+r} \cg_-^{j} \amalg \cb $.

Finally, we consider ~$so^*(2n)$. First we suppose ~$n=2r$. Then ~$\cm_0 ~=~   so(3) \oplus \cdots \oplus so(3)$, ($r$ factors), ~$\dim_\bbr \ca_0 ~=~r$, ~$\dim_\bbr \cn_0 ~=~r(4r-3)$.
Note also that
~$\cm_{0h}^\bbc \oplus \ca_0^\bbc \cong \ch(so(2n,\bbc)$.\\
Thus, we have ~$\cp^\bbc_0 \cong P_{1,3,\ldots,n-1} \cong \amalg_{i=1}^r \cg_-^{2i-1}\amalg \cb $.\\
Next we suppose ~$n=2r+1$. Then ~$\cm_0 ~=~  so(2)\oplus so(3) \oplus \cdots \oplus so(3)$, ($r$ factors),
~$\dim_\bbr \ca_0 ~=~r$, ~$\dim_\bbr \cn_0 ~=~r(4r+1)$.
Thus, we have ~$\cp^\bbc_0 \cong P_{1,3,\ldots,n-2} \cong \amalg_{i=1}^r \cg_-^{2i-1}\amalg \cb $

For the lack of space we leave consideration of the  exceptional  real Lie algebras for a subsequent publication
 \cite{VKD-prep}.

\noindent {\bf Discussion:}~ Another main ingredient of our approach as follows. We group
the (reducible) ERs with the same Casimirs in sets called ~{\it
multiplets} \cite{Dob}. The multiplet corresponding to fixed values of the
Casimirs may be depicted as a connected graph, the {\it vertices} of
which correspond to the reducible ERs and the {\it lines (arrows)}
between the vertices correspond to intertwining operators.  The
explicit parametrization of the multiplets and of their ERs is
important for understanding of the situation. The notion of multiplets was introduced in \cite{Dobmul}. %
   Then it was applied to  (infinite-dimensional) (super-)algebras,
quantum groups and other symmetry objects. For a current summary of
these developments, see \cite{VKD1,VKD2}, for further developments - \cite{VKD-prep}.

\medskip

\noindent {\bf Acknowledgment}.  The author has received partial support from  Bulgarian NSF Grant   DN-18/1.

\end{document}